\documentclass[11pt]{article}

\usepackage{latexsym}
\usepackage{amsfonts}
\usepackage{amssymb}
\usepackage{mathrsfs}
\usepackage[all]{xy}
\usepackage{epic}
\usepackage{eepic}

\newcommand{\blda}{\mathbf{A}}

\newcommand{\bldz}{\mathbf{z}}

\newcommand{\sgn}{\mbox{sign\,}}

\newcommand{\be}{\begin{equation}}
\newcommand{\ee}{\end{equation}}
\newcommand{\bea}{\begin{eqnarray}}
\newcommand{\eea}{\end{eqnarray}}
\newcommand{\bean}{\begin{eqnarray*}}
\newcommand{\eean}{\end{eqnarray*}}
\newcommand{\brray}{\begin{array}}
\newcommand{\erray}{\end{array}}

\newcommand{\newsection}[1]{\setcounter{equation}{0}
\setcounter{dfn}{0}
\section{#1}}

\newtheorem{dfn}{Definition}[section]
\newtheorem{thm}[dfn]{Theorem}
\newtheorem{lmma}[dfn]{Lemma}
\newtheorem{ppsn}[dfn]{Proposition}
\newtheorem{crlre}[dfn]{Corollary}
\newtheorem{xmpl}[dfn]{Example}
\newtheorem{rmrk}[dfn]{Remark}

\newcommand{\bdfn}{\begin{dfn}\rm}
\newcommand{\bthm}{\begin{thm}}
\newcommand{\blmma}{\begin{lmma}}
\newcommand{\bppsn}{\begin{ppsn}}
\newcommand{\bcrlre}{\begin{crlre}}
\newcommand{\bxmpl}{\begin{xmpl}}
\newcommand{\brmrk}{\begin{rmrk}\rm}

\newcommand{\edfn}{\end{dfn}}
\newcommand{\ethm}{\end{thm}}
\newcommand{\elmma}{\end{lmma}}
\newcommand{\eppsn}{\end{ppsn}}
\newcommand{\ecrlre}{\end{crlre}}
\newcommand{\exmpl}{\end{xmpl}}
\newcommand{\ermrk}{\end{rmrk}}

\newcommand{\bbz}{\mathbb{Z}}
\newcommand{\bbn}{\mathbb{N}}

\newcommand{\fks}{\mathfrak{S}}
\newcommand{\scrc}{\mathscr{C}}
\newcommand{\scrf}{\mathscr{F}}

\newcommand{\cla}{\mathcal{A}}

\newcommand{\clh}{\mathcal{H}}

\newcommand{\clk}{\mathcal{K}}

\newcommand{\clg}{\mathcal{G}}

\newcommand{\rchi}{\raisebox{.4ex}{\ensuremath{\chi}}}

\newcommand{\prf}{\noindent{\it Proof\/}: }

\newcommand{\seq}{\subseteq}

\newcommand{\dmove}{\textsl{move}\ }

\def \qed { \mbox{}\hfill
$\Box$\vspace{1ex}}

\newcommand{\half}{\frac{1}{2}}

\begin{document}


\author{{\sc Partha Sarathi Chakraborty} and
{\sc Arupkumar Pal}}
\title{Characterization of spectral triples: A combinatorial approach}
\maketitle
 \begin{abstract}
   We describe a general technique to study Dirac operators on
   noncommutative spaces under some additional assumptions.  The main
   idea is to capture the compact resolvent condition in a
   combinatorial set up.  Using this, we then prove that for
   a certain class of representations of
   the $C^*$-algebra $C(SU_q(\ell+1))$, any Dirac operator that
   diagonalises with respect to the natural basis of the underlying
   Hilbert space must have trivial sign.
 \end{abstract}
{\bf AMS Subject Classification No.:} {\large 58}B{\large 34}, {\large
46}L{\large 87}, {\large
  19}K{\large 33}



\newsection{Introduction}
A spectral triple is the starting
point in noncommutative geometry (NCG) where a geometric space is
described by a triple $(\cla,\clh,D)$, with $\cla$ being an involutive
algebra represented as bounded operators on a Hilbert space $\clh$,
and $D$ being a selfadjoint operator with compact resolvent and having
bounded commutators with the algebra elements. This $D$ should be
nontrivial in the sense that the associated Kasparov module should
give a nontrivial element in $K$-homology.
Observe that the self-adjoint operator $D$ in a spectral
triple comes with two very crucial restrictions on it, namely, it has to have
compact resolvent, and must have bounded commutators with algebra
elements. Various analytic consequences of the compact resolvent
condition (growth properties of the commutators of the algebra
elements with the sign of $D$) have been used in the past by various
authors.  Here we will take a new approach that will help us exploit
it from a combinatorial point of view.  The idea is very simple. Given a
selfadjoint operator with compact resolvent, one can associate with it
a certain  graph in a natural way.  This makes it possible to do a detailed
combinatorial analysis of the growth restrictions (on the eigenvalues
of $D$) that come from the boundedness of the commutators, and to
characterize the sign of the operator $D$ completely.

In the next section, we will outline the strategy.
It should be noted that this technique has already
been used implicitly in characterizing spectral
triples for the quantum $SU(2)$ group by the
authors in~\cite{c-p1} and~\cite{c-p2}.
Here we will present the scheme in a more explicit
way and use it  in the remaining sections
to study a more complicated and important case.
The case that we treat is analogous to the one for $SU_q(2)$
treated in~\cite{c-p2}.
We will take a large class of representations of the
 $C^*$-algebra $C(SU_q(\ell+1))$, which
includes the irreducibles in particular, and use the
general scheme described in section~1 to prove that
for a large majority of
these representations, any Dirac operator that
diagonalises nicely with respect to
the canonical orthonormal basis must have trivial
sign.


\newsection{The combinatorial set up}\label{g_scheme}
Suppose
$\cla$ is a $C^*$-algebra represented on a Hilbert space, and suppose
we want to have an idea about all operators $D$ that will make
$(\cla,\clh,D)$ into a spectral triple.
Of course, in this
generality, the problem would be intractable in most cases. We will
impose some extra conditions on this $D$ that will be natural
from the context.
This would give some information about the spectral resolution
$\sum_{\gamma\in\Gamma}d_\gamma P_\gamma$, more
specifically some idea
about the set $\Gamma$ and
a diagonalising basis for
$D$. Note that since $D$ is known to be self-adjoint with discrete
spectrum, there always exists such a basis.
Next, let $c$ be a positive real. Construct a graph $\clg_c$
by taking the vertex set $V$ to be $\Gamma$ and  by
joining two points $\gamma$ and
$\gamma'$ in $V=\Gamma$ by an edge if
$|d_\gamma-d_{\gamma'}|<c$.
Define
$V^+=\{\gamma\in V: d_\gamma>0\}$ and
$V^-=\{\gamma\in V: d_\gamma<0\}$.
 One can assume
without loss in generality that the null space of $D$ is trivial,
as this can be achieved just by a compact perturbation.
Thus $(V^+,V^-)$ gives us  a partition of the vertex set
$\Gamma=V$.
Call  two paths $(v_1,v_2,\ldots,v_m)$ and $(w_1,w_2,\ldots,w_n)$
in $\clg_c$
(or more generally in any graph $\clg$)
\textbf{disjoint} if
the sets $\{v_1,v_2,\ldots,v_m\}$ and
$\{w_1,w_2,\ldots,w_n\}$ do not intersect.
Now observe that there can not exist
infinitely many disjoint paths from $V^+$ to $V^-$.
This is because if $(v_1,v_2,\ldots,v_m)$ is a
path from $V^+$ to $V^-$, then for some $v_i$, one
must have $d_{v_i}\in[-c,c]$. Thus if there are
infinitely many disjoint paths from $V^+$ to $V^-$,
it would contradict the compact resolvent condition
on $D$.

We say that a partition $(V_1,V_2)$ of the vertex set
in a graph admits an \textbf{infinte ladder} if there are
infinitely many disjoint paths from $V_1$ to $V_2$.
We call a partition $(V_1,V_2)$ \textbf{sign-determining} if
it does not admit an infinite ladder.
Thus the partition $(V^+,V^-)$ of the vertex set
in the graph $\clg_c$ is sign-determining.

Since we do not know the operator $D$ apriori,
we proceed from the other direction.
Using the action of the algebra elements on the basis
elements of $\clh$ and using the boundedness of their
commutators with $D$, we get certain growth restrictions
on the $d_\gamma$'s. These will give us some information
about the edges in the graph. We exploit this knowledge
to characterize those partitions $(V_1,V_2)$ of the vertex set
that do not admit any infinite ladder.
This amounts to characterizing the sign of the operator $D$,
in the sense that the sign of $D$ must be of the form
$\sum_{\gamma\in V_1}P_\gamma-\sum_{\gamma\in V_2}P_\gamma$.

Note here that whether or not a partition admits
an infinite ladder will depend on the value of $c$.
For a specific value of $c$, the graph
$\clg_c$ may have no edges, or too few edges (if the singular values
of $D$ happen to grow too fast). In such a case, there would
exist too many partitions that are sign-determining, and as a result,
will not be very useful. Therefore we will be interested only
in those partitions that remain sign-determining
for all sufficiently large $c$.



\newsection{The group $SU_q(\ell+1)$}
Let $\mathfrak{g}$ be a complex simple Lie algebra of rank $\ell$.
let $(\!(a_{ij})\!)$ be the associated Cartan matrix,
$q$ be a real number lying in the interval $(0,1)$
and let $q_i=q^{(\alpha_i,\alpha_i)/2}$, where $\alpha_i$'s are the simple roots
of $\mathfrak{g}$.
Then the quantised universal envelopping algebra (QUEA)
$U_q(\mathfrak{g})$ is the algebra
generated by $E_i$, $F_i$, $K_i$ and $K_i^{-1}$, $i=1,\ldots,\ell$, satisfying the
following relations
\begin{displaymath}
K_iK_j=K_jK_i,\quad K_iK_i^{-1}=K_i^{-1}K_i=1,
\end{displaymath}
\begin{displaymath}
K_iE_jK_i^{-1}=q_i^{\half a_{ij}}E_j,\quad
K_iF_jK_i^{-1}=q_i^{-\half a_{ij}}F_j,
\end{displaymath}
\begin{displaymath}
E_iF_j-F_jE_i=\delta_{ij}\frac{K_i^2-K_i^{-2}}{q_i-q_i^{-1}},
\end{displaymath}
\begin{displaymath}
\sum_{r=0}^{1-a_{ij}}(-1)^r{{1-a_{ij}}\choose r}_{q_i}
  E_i^{1-a_{ij}-r}E_jE_i^r =0 \quad\forall\, i\neq j,
\end{displaymath}
\begin{displaymath}
\sum_{r=0}^{1-a_{ij}}(-1)^r{{1-a_{ij}}\choose r}_{q_i}
  F_i^{1-a_{ij}-r}F_jF_i^r =0\quad \forall\, i\neq j,
\end{displaymath}
where ${n\choose r}_q$ denote the $q$-binomial coefficients.
Hopf *-structure comes from the following maps:
\[
\Delta(K_i)=K_i\otimes K_i,\quad \Delta(K_i^{-1})=K_i^{-1}\otimes K_i^{-1},
\]
\[
\Delta(E_i)=E_i\otimes K_i  + K_i^{-1}\otimes E_i,\quad
\Delta(F_i)=F_i\otimes K_i  + K_i^{-1}\otimes F_i,
\]
\[
\epsilon(K_i)=1,\quad \epsilon(E_i)=0=\epsilon(F_i),
\]
\[
S((K_i)=K_i^{-1},\quad S(E_i)=-q_iE_i,\quad S(F_i)=-q_i^{-1}F_i,
\]
\[
K_i^*=K_i,\quad E_i^*=-q_i^{-1}F_i,\quad F_i^*=-q_iE_i.
\]

In the type A case, the associated Cartan matrix is given by
\[
a_{ij}=\cases{2& if $i=j$,\cr
              -1 & if $i=j\pm1$,\cr
               0 & otherwise,}
\]
and $(\alpha_i,\alpha_i)=2$ so that $q_i=q$ for all $i$.
The QUEA in this case is denoted by $u_q(su(\ell+1))$.

Take the collection of matrix entries of all finite-dimensional
unitarizable $u_q(su(\ell+1))$-modules. The algebra generated by these
gets a natural Hopf*-structure as the dual of $u_q(su(\ell+1))$.  One
can also put a natural $C^*$-norm on this.  Upon completion with
respect to this norm, one gets a unital $C^*$-algebra that plays the
role of the algebra of continuous functions on $SU_q(\ell+1)$.  For a
detailed account of this, refer to chapter~3, \cite{ko-so}.  In
\cite{w}, Woronowicz gave a different description of this
$C^*$-algebra.  which was later shown by Rosso (\cite{r}) to be
equivalent to the earlier one.

For remainder of this article, we will take $G$ to be $SU_q(\ell+1)$
and $\cla$ will be the $C^*$-algebra of continuous functions on $G$.

\newsection{Irreducible representations}
 All irreducible representations of
the $C^*$-algebra $\cla$ are well-known (\cite{ko-so}).
Let us briefly recall those here.
The Weyl group for $SU_q(\ell+1)$ is isomorphic to the
permutations group $\fks_{\ell+1}$ on $\ell+1$ symbols.
Denote by $s_i$ the transposition $(i,i+1)$.
Then $\{s_1,\ldots,s_\ell\}$ form a set of generators
for $\fks_{\ell+1}$. Any $\omega\in\fks_{\ell+1}$
can be written as a product
\[
\omega=(s_{k_\ell}s_{k_\ell+1}\ldots s_\ell)
 (s_{k_{\ell-1}}s_{k_{\ell-1}+1}\ldots s_{\ell-1})
\ldots (s_{k_2}s_2)(s_{k_1}),
\]
where $k_i$'s are integers satisfying
$0\leq k_i\leq i$, with the understanding that
$k_i=0$ means that the string
$(s_{k_i}s_{k_i+1}\ldots s_i)$
is missing.
It follows from the strong deletion condition
in the characterization of Coxeter system by Tits
(see~\cite{g}) that
the expression for $\omega$ given above is a reduced
word in the generators $s_i$. We will denote the length
of an element $\omega$ by $\ell(\omega)$.

Let $S$ and $N$ be the following operators on $L_2(\bbz)$:
\[
Se_n=e_{n-1},\quad Ne_n=ne_n.
\]
We will denote by the same symbols their restrictions to
$L_2(\bbn)$ whenever there is no chance of ambiguity.
Denote by $\psi_{s_i}$ the following representation of $\cla$
on $L_2(\bbn)$:
\[
\psi_{s_i}(u_{ab})=\cases{
                   \sqrt{I-q^{2N+2}}S & if $a=b=i$,\cr
                  S^*\sqrt{I-q^{2N+2}} &if $a=b=i+1$,\cr
                   -q^{N+1} & if $a=i$, $b=i+1$,\cr
                    q^N   & if $a=i+1$, $b=i$,\cr
                    \delta_{ab}I & otherwise.}
\]
Now suppose $\omega\in\fks_{\ell+1}$ is given by
$s_{i_1}s_{i_2}\ldots s_{i_k}$.
Define $\psi_\omega$ to be
$\psi_{s_{i_1}}\ast\psi_{s_{i_2}}\ast\ldots\ast\psi_{s_{i_k}}$
(for two representations $\phi$ and $\psi$, $\phi\ast\psi$ denote the representation
$(\phi\otimes\psi)\Delta$).

Next, let $\bldz=(z_1,\ldots,z_\ell)\in (S^1)^\ell$.
Define
\[
\rchi_\bldz(u_{ab})=\cases{z_a\delta_{ab} & if $a=1$,\cr
                         \bar{z}_\ell\delta_{ab} & if $a=\ell+1$,\cr
                         \bar{z}_{a-1} z_a \delta_{ab} & otherwise.}
\]
Define $\rchi$ to be the integral $\int_{\bldz\in (S^1)^\ell}\rchi_\bldz d\bldz$.
Finally, define $\pi_{\omega,\bldz}=\psi_\omega\ast\rchi_\bldz$ and
$\pi_\omega=\psi_\omega\ast\rchi$.
It is known (\cite{ko-so}) that $\pi_{\omega,\bldz}$'s constitute
all the ireducible representations of the $C^*$-algebra $\cla$.

Let us introduce a few notations that will be handy later.
For a subset $\Lambda=\{i_1,\ldots,i_k\}\seq\{1,2,\ldots,\ell\}$,
where $i_1<i_2<\ldots < i_k$, denote by $s_\Lambda$ the element
$s_{i_1}s_{i_2}\ldots s_{i_k}$ of $\fks_{\ell+1}$.
Call a subset $J$ of $\{1,2,\ldots,\ell\}$ an
\textbf{interval} if it is of the form
$\{j,j+1,\ldots,j+s\}$.
Then for any element $\omega$ of the Weyl group,
there are intervals $\Lambda_1, \Lambda_2,\ldots,\Lambda_t$
with $\max \Lambda_r>\max \Lambda_s$ for $r>s$ such that
\be
\omega=s_{\Lambda_t}s_{\Lambda_{t-1}}\ldots s_{\Lambda_1}.
\ee
Moreover, as long as we demand that $\Lambda_j$'s are intervals
and obey $\max \Lambda_r>\max \Lambda_s$ for $r>s$, an element $\omega$
determines the subsets $\Lambda_t, \ldots,\Lambda_1$ uniquely.
Let $\Lambda$ be the disjoint union of the $\Lambda_j$'s,
that is, $\Lambda=\cup_{j=1}^t\{(j,i): i\in \Lambda_j\}$.
Write  $\Lambda_0=\{(0,i): i=1,2,\ldots,\ell\}$.
Often we will identify $\Lambda_0$ with the set $\{1,2,\ldots,\ell\}$.
Let $\Gamma=\bbn^\Lambda\times \bbz^{\Lambda_0}
= \bbn^{\Lambda_t}\times \bbn^{\Lambda_{t-1}}\times\ldots\times \bbn^{\Lambda_1}\times \bbz^{\Lambda_0}$.
The Hilbert space on which $\pi_\omega$ acts is $L_2(\Gamma)$.
We will denote by $\{e_\gamma:\gamma\in\Gamma\}$ the
canonical orthonormal basis for this Hilbert space.

\newsection{Diagram representation of $\pi_\omega$}
Let us describe how to use a diagram to represent the irreducible $\psi_{s_i}$.
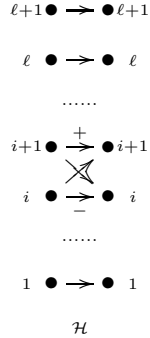
\begin{center}
\begin{minipage}{200pt}
  \def\labelstyle{\scriptscriptstyle}
  \xymatrix@C=10pt@R=8pt{\ar@{}[r]|{\ell+1}&{\bullet}\ar[r]
&{\bullet}\ar@{}[r]|{\ell+1}& \\
    \ar@{}[r]|\ell&{\bullet}\ar[r] &{\bullet}\ar@{}[r]|{\ell}& \\
    &\ar@{}[r]|{\ldots\ldots}&  \\
    \ar@{}[r]|{i+1}&{\bullet}\ar@{->}[r]^{+}\ar@{->}[dr] &{\bullet}\ar@{}[r]|{i+1}&
\\
    \ar@{}[r]|i&{\bullet}\ar@{->}[r]_{-}\ar@{->}[ur] &{\bullet}\ar@{}[r]|{i}& \\
    &\ar@{}[r]|{\ldots\ldots}&  \\
    \ar@{}[r]|1&{\bullet}\ar[r] &{\bullet}\ar@{}[r]|{1}& \\
    &\ar@{}[r]|{\clh} && }
\end{minipage}
     \end{center}
In this diagram, each path from a node $k$ on the
left to a node $l$ on the right stands for an
operator on $\clh=L_2(\bbn)$. A horizontal
unlabelled line stands for the identity operator,
a horizontal line labelled with a $+$ sign
       stands for $S^*\sqrt{I-q^{2N+2}}$ and
one labelled with a $-$ sign stands for
$\sqrt{I-q^{2N+2}}S$. A diagonal
line going upward represents $-q^{N+1}$ and a diagonal line going
downward represents $q^N$.
Now $\psi_{s_i}(u_{kl})$ is the operator represented by the path from $k$ to $l$,
and is zero if there is no such path.
Thus, for example, $\psi_{s_i}(u_{11})$
is $I$, $\psi_{s_i}(u_{12})$ is zero, whereas $\psi_{s_i}(u_{ii+1})=-q^{N+1}$, if $i>1$.

Next, let us explain how to represent $\psi_{s_i}\ast \psi_{s_j}$.
Simply put the two diagrams representing $\psi_{s_i}$ and $\psi_{s_j}$
adjacent to each other, and identify, for each row, the node on the right side
of the diagram for $\psi_{s_i}$ with the node on the left in the diagram
for $\psi_{s_j}$. Now, $\psi_{s_i}\ast \psi_{s_j}(u_{kl})$ would be
an operator on $L_2(\bbn)\otimes L_2(\bbn)$ determined by all the paths from
the node $k$ on the left to the node $l$ on the right. It would be zero if
there is no such path and if there are more than one paths, then it would be the sum of
the operators given by each such path. Thus, we have the following operation
on the elementary diagrams described above:
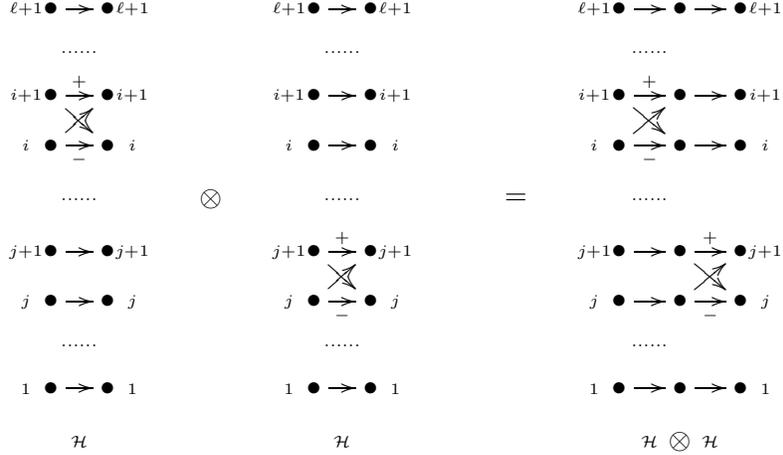
\begin{center}
\begin{minipage}{400pt}
  \def\labelstyle{\scriptscriptstyle}
  \xymatrix@C=10pt@R=8pt{
    \ar@{}[r]|{\ell+1}&{\bullet}\ar[r] &{\bullet}\ar@{}[r]|{\ell+1}&&&
         \ar@{}[r]|{\ell+1}&{\bullet}\ar[r] &{\bullet}\ar@{}[r]|{\ell+1}&&&&
  \ar@{}[r]|{\ell+1}&{\bullet}\ar[r] &{\bullet}\ar[r] &{\bullet}\ar@{}[r]|{\ell+1}&\\
    &\ar@{}[r]|{\ldots\ldots}&&&&&
       \ar@{}[r]|{\ldots\ldots}&&&&&&
           \ar@{}[r]|{\ldots\ldots}&  \\
    \ar@{}[r]|{i+1}&{\bullet}\ar@{->}[r]^{+}\ar@{->}[dr]&{\bullet}\ar@{}[r]|{i+1}
&&&
     \ar@{}[r]|{i+1}&{\bullet}\ar[r] &{\bullet}\ar@{}[r]|{i+1}&&&&
      \ar@{}[r]|{i+1}&{\bullet}\ar@{->}[r]^{+}\ar@{->}[dr]&{\bullet}\ar[r]
&{\bullet}\ar@{}[r]|{i+1}& \\
    \ar@{}[r]|i&{\bullet}\ar@{->}[r]_{-}\ar@{->}[ur]&{\bullet}\ar@{}[r]|{i} &&&
       \ar@{}[r]|i&{\bullet}\ar[r] &{\bullet}\ar@{}[r]|{i}&&&&
     \ar@{}[r]|i&{\bullet}\ar@{->}[r]_{-}\ar@{->}[ur]&{\bullet}\ar[r]
&{\bullet}\ar@{}[r]|{i}& \\
    &\ar@{}[r]|{\ldots\ldots}&&&{\otimes}&&
        \ar@{}[r]|{\ldots\ldots}&&&&{=}&&
           \ar@{}[r]|{\ldots\ldots}&  \\
    \ar@{}[r]|{j+1}&{\bullet}\ar[r]&{\bullet}\ar@{}[r]|{j+1} &&&
    \ar@{}[r]|{j+1}&{\bullet}\ar@{->}[r]^{+}\ar@{->}[dr]
&{\bullet}\ar@{}[r]|{j+1}&&&&
 \ar@{}[r]|{j+1}&{\bullet}\ar[r]&{\bullet}\ar@{->}[r]^{+}\ar@{->}[dr]
&{\bullet}\ar@{}[r]|{j+1}& \\
    \ar@{}[r]|j&{\bullet}\ar[r]&{\bullet}\ar@{}[r]|{j} &&&
 \ar@{}[r]|j&{\bullet}\ar@{->}[r]_{-}\ar@{->}[ur] &{\bullet}\ar@{}[r]|{j}&&&&
 \ar@{}[r]|j&{\bullet}\ar[r]
&{\bullet}\ar@{->}[r]_{-}\ar@{->}[ur]&{\bullet}\ar@{}[r]|{j}& \\
    &\ar@{}[r]|{\ldots\ldots}&&&&&
      \ar@{}[r]|{\ldots\ldots}&&&&&&
         \ar@{}[r]|{\ldots\ldots}&  \\
    \ar@{}[r]|1&{\bullet}\ar[r]&{\bullet}\ar@{}[r]|{1} &&&
     \ar@{}[r]|1&{\bullet}\ar[r] &{\bullet}\ar@{}[r]|{1}&&&&
    \ar@{}[r]|1&{\bullet}\ar[r]&{\bullet}\ar[r] &{\bullet}\ar@{}[r]|{1}&  \\
    &\ar@{}[r]|{\clh}& &&&
     &\ar@{}[r]|{\clh} &&&&&
    &\ar@{}[r]|{\clh}&{\otimes}\ar@{}[r]|{\clh} &&
    }
\end{minipage}
     \end{center}

Next, we come to $\rchi$. The underlying Hilbert space now is
$L_2(\bbz^{\Lambda_{0}})\cong L_2(\bbz)^{\otimes \ell}$ (to avoid any ambiguity, we have used
hollow circles to denote the nodes as opposed to the bullets
used in the earlier case); an unlabelled horizontal arrow stands for
$I$ in the corresponding component of $L_2(\bbz)^{\otimes \ell}$,
an arrow labelled with a `$+$' above it indicates $S^*$ and one
labelled `$-$' below it stands for $S$.
As earlier, $\rchi(u_{kl})$ stands for the operator
on $L_2(\bbz)^{\otimes \ell}$ represented by the path
from $k$ on the left to $l$ on the right.
In the diagram below, $\clk$ will stand for $L_2(\bbz)$.
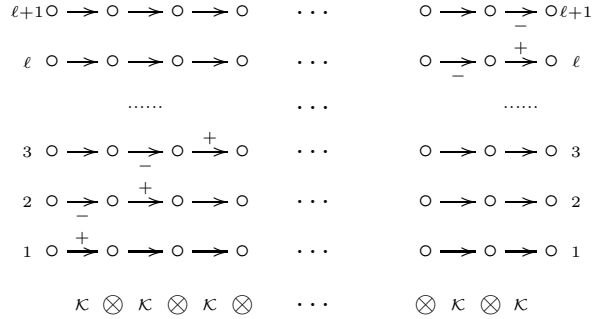
\begin{center}
\begin{minipage}{200pt}
  \def\labelstyle{\scriptscriptstyle}
 \xymatrix@C=10pt@R=8pt{
    \ar@{}[r]|{\ell+1}&{\circ}\ar[r]&{\circ}\ar[r]&{\circ}\ar[r]&{\circ}&{\ldots}
       &&{\circ}\ar[r]&{\circ}\ar@{->}[r]_{-}&{\circ}\ar@{}[r]|{\ell+1}&  \\
    \ar@{}[r]|\ell&{\circ}\ar[r]&{\circ}\ar[r]&{\circ}\ar[r]&{\circ}&{\ldots}
       &&{\circ}\ar@{->}[r]_{-}&{\circ}\ar@{->}[r]^{+}&{\circ}\ar@{}[r]|{\ell}& \\
    &&\ar@{}[r]|{\ldots\ldots}&&&{\ldots}&&&\ar@{}[r]|{\ldots\ldots}&&  \\
    \ar@{}[r]|3&{\circ}\ar[r]&{\circ}\ar@{->}[r]_{-}&{\circ}\ar@{->}[r]^{+}&{\circ}&{\ldots}
       &&{\circ}\ar[r]&{\circ}\ar[r]&{\circ}\ar@{}[r]|{3}& \\
    \ar@{}[r]|2&{\circ}\ar@{->}[r]_{-}&{\circ}\ar@{->}[r]^{+}&{\circ}\ar[r]&{\circ}&{\ldots}
       &&{\circ}\ar[r]&{\circ}\ar[r]&{\circ}\ar@{}[r]|{2}& \\
    \ar@{}[r]|1&{\circ}\ar@{->}[r]^{+}&{\circ}\ar[r]&{\circ}\ar[r]&{\circ}&{\ldots}
       &&{\circ}\ar[r]&{\circ}\ar[r]&{\circ}\ar@{}[r]|{1}& \\
    &\ar@{}[r]|{\clk}&{\otimes}\ar@{}[r]|{\clk}&{\otimes}\ar@{}[r]|{\clk}&{\otimes}&{\ldots}
       &&{\otimes}\ar@{}[r]|{\clk}&{\otimes}\ar@{}[r]|{\clk}&&
}
\end{minipage}
     \end{center}

Finally, we come to the description of $\pi_\omega$.
As we have already remarked, reduced expression for $\omega$
is of the form
$\omega=(s_{k_n}s_{k_n+1}\ldots s_n)(s_{k_{n-1}}\ldots s_{n-1})\ldots (s_{k_2}s_2)(s_{k_1})$.
To get the diagram for $\pi_\omega$, we simply put
the diagram for
$\psi_{s_{k_n}}\ast\ldots\ast\psi_{s_{k_1}}$
and that for $\rchi$ side by side and identify the nodes
on the right of the first diagram with the corresponding ones on the left
of the second diagram. Thus for example, if
 $\omega=(s_2s_3s_4)(s_3)(s_1s_2)(s_1)$, then the following
diagram represents $\pi_\omega$:


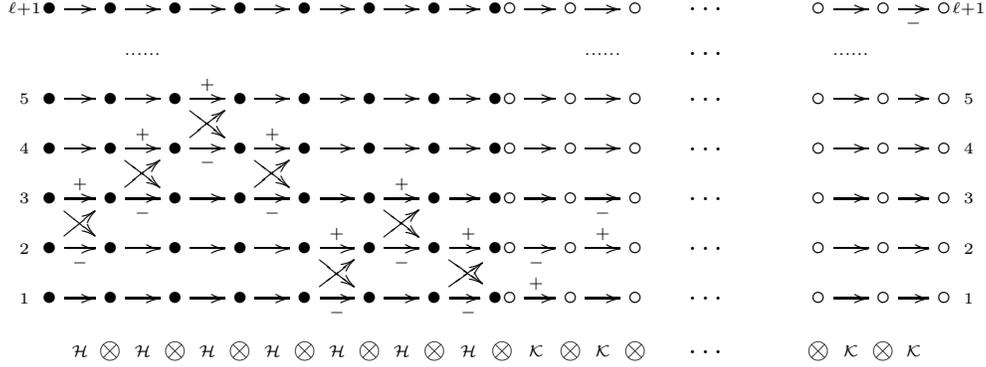
\begin{center}
\begin{minipage}{450pt}
  \def\labelstyle{\scriptscriptstyle}
  \xymatrix@C=10pt@R=8pt{
    \ar@{}[r]|{\ell+1}&{\bullet}\ar[r]&{\bullet}\ar[r]&{\bullet}\ar[r]&
          {\bullet}\ar[r]&{\bullet}\ar[r]&{\bullet}\ar[r]&{\bullet}\ar[r]&
{\bullet}{\circ}\ar[r]&{\circ}\ar[r]&{\circ}&{\ldots}
       &&{\circ}\ar[r]&{\circ}\ar@{->}[r]_{-}&{\circ}\ar@{}[r]|{\ell+1}&  \\
   &&\ar@{}[r]|{\ldots\ldots}&&&&& &&\ar@{}[r]|{\ldots\ldots}&&
         {\ldots}&&\ar@{}[r]|{\ldots\ldots}&&&  \\
    \ar@{}[r]|5&{\bullet}\ar[r]&{\bullet}\ar[r]&{\bullet}\ar@{->}[r]^{+}\ar[dr]&
          {\bullet}\ar[r]&{\bullet}\ar[r]&{\bullet}\ar[r]&{\bullet}\ar[r]&
{\bullet}{\circ}\ar[r]&{\circ}\ar[r]&{\circ}&{\ldots}
       &&{\circ}\ar[r]&{\circ}\ar[r]&{\circ}\ar@{}[r]|{5}& \\
    \ar@{}[r]|4&{\bullet}\ar[r]&{\bullet}\ar@{->}[r]^{+}\ar[dr]&{\bullet}\ar@{->}[r]_{-}\ar[ur]&
          {\bullet}\ar@{->}[r]^{+}\ar[dr]&{\bullet}\ar[r]&{\bullet}\ar[r]&{\bullet}\ar[r]&
{\bullet}{\circ}\ar[r]&{\circ}\ar[r]&{\circ}&{\ldots}
       &&{\circ}\ar[r]&{\circ}\ar[r]&{\circ}\ar@{}[r]|{4}& \\
    \ar@{}[r]|3&{\bullet}\ar@{->}[r]^{+}\ar[dr]&{\bullet}\ar@{->}[r]_{-}\ar[ur]&{\bullet}\ar[r]&
          {\bullet}\ar@{->}[r]_{-}\ar[ur]&{\bullet}\ar[r]&{\bullet}\ar@{->}[r]^{+}\ar[dr]&
     {\bullet}\ar[r]&
{\bullet}{\circ}\ar[r]&{\circ}\ar@{->}[r]_{-}&{\circ}&{\ldots}
       &&{\circ}\ar[r]&{\circ}\ar[r]&{\circ}\ar@{}[r]|{3}& \\
    \ar@{}[r]|2&{\bullet}\ar@{->}[r]_{-}\ar[ur]&{\bullet}\ar[r]&{\bullet}\ar[r]&
{\bullet}\ar[r]&{\bullet}\ar@{->}[r]^{+}\ar[dr]&{\bullet}\ar@{->}[r]_{-}\ar[ur]&
    {\bullet}\ar@{->}[r]^{+}\ar[dr]&
{\bullet}{\circ}\ar@{->}[r]_{-}&{\circ}\ar@{->}[r]^{+}&{\circ}&{\ldots}
       &&{\circ}\ar[r]&{\circ}\ar[r]&{\circ}\ar@{}[r]|{2}& \\
    \ar@{}[r]|1&{\bullet}\ar[r]&{\bullet}\ar[r]&{\bullet}\ar[r]&
          {\bullet}\ar[r]&{\bullet}\ar@{->}[r]_{-}\ar[ur]&{\bullet}\ar[r]&
        {\bullet}\ar@{->}[r]_{-}\ar[ur]&
{\bullet}{\circ}\ar@{->}[r]^{+}&{\circ}\ar[r]&{\circ}&{\ldots}
       &&{\circ}\ar[r]&{\circ}\ar[r]&{\circ}\ar@{}[r]|{1}& \\
    &\ar@{}[r]|{\clh}&{\otimes}\ar@{}[r]|{\clh}&{\otimes}\ar@{}[r]|{\clh}&
          {\otimes}\ar@{}[r]|{\clh}&{\otimes}\ar@{}[r]|{\clh}&{\otimes}\ar@{}[r]|{\clh}&
        {\otimes}\ar@{}[r]|{\clh}&
{\otimes}\ar@{}[r]|{\clk}&{\otimes}\ar@{}[r]|{\clk}&{\otimes}&{\ldots}
       &&{\otimes}\ar@{}[r]|{\clk}&{\otimes}\ar@{}[r]|{\clk}&&
  }
\end{minipage}
     \end{center}

The diagram for $\pi_\omega$ introduced above will
play an important role in what follows.


\newsection{Boundedness of commutators}
Our goal is to study operators $D$ on the space $\clh_\omega=L_2(\Gamma)$
that diagonalize with respect to the natural canonical basis,
and makes $(\pi_\omega(\cla),\clh_\omega, D)$ a spectral triple.
Since $D$ is a self-adjoint operator with discrete spectrum,
it is of the form
$\sum_{\gamma\in\Gamma}d(\gamma)e_\gamma$,
where $d(\gamma)\geq 0$ for all $\gamma$.

\bdfn\label{neqmove}
A \textbf{move} will mean a path from a node on the left to a node
on the right in the diagram representing $\pi_\omega$.
More formally, a \dmove is a $(t+1)$-tuple of pairs
$((i_t,j_t),\ldots,(i_0,j_0))$ such that
\begin{enumerate}
\item $j_k=i_{k-1}$ for $k\geq 1$, $i_0=j_0$,
\item for $k\geq 1$, $j_k<i_k$ implies $j_k=i_k-1$ and $j_k\in \Lambda_k$,
\item for $k\geq 1$, $j_k>i_k$ implies $i_k,i_k+1,\ldots,j_k-1\in \Lambda_k$.
\end{enumerate}
(the pair $(i_k,j_k)$ will be referred as the $k$\raisebox{.4ex}{th}
segment of the \dmove lying in the $k$\raisebox{.4ex}{th} string from the right).\\
Observe that the above three conditions imply
in particular that $1\leq i_k,j_k\leq \ell+1$
for all $k$.
We will use the special notation $H_r$ for the \dmove for which
each $i_k$ and $j_k$ equals $r$.
\edfn

Given a \dmove $p$, we will next define an element
$m_p\in\bbz^\Lambda\times\bbz^{\Lambda_0}=
  \bbz^{\Lambda\cup\Lambda_0}$
whose coordinates are all 0 or $\pm 1$.
Let $p=((i_t,j_t),\ldots,(i_0,j_0))$.
Define $m_p$ by the following prescription:
\[
m_p(r,s)=\cases{ -1 & if $r=0, s=i_0-1$ or $r\geq 1, s=j_r\geq i_r$,\cr
    +1 & if $r=0, s=i_0$ or $r\geq 1, j_r\geq i_r=s+1$,\cr
                 0 & otherwise.}
\]
Thus for $r>0$, $m_p(r,\cdot)$ will look like
\[
\cases{(0,0,\ldots,0) & if $i_r<\min \Lambda_r$ or $i_r>\max \Lambda_r+1$ or $j_r=i_r-1$,\cr
(-1,0,\ldots,0) & if $i_r=j_r=\min \Lambda_r$,\cr
(0,0,\ldots,0,1) & if $i_r=j_r=\max \Lambda_r+1$,\cr
(\underbrace{0,\ldots,0}_{i_r-2},1,
   \underbrace{0,\ldots,0}_{j_r-i_r},-1,0,\ldots,0) &
                        if $\min \Lambda_r<i_r\leq j_r$,}
\]
and $m_p(0,\cdot)$ will be of the form
\[
\cases{(1,0,\ldots,0) & if $i_0=j_0=1$,\cr
   (0,0,\ldots,0,-1) & if $i_0=j_0=\ell+1$,\cr
   (\underbrace{0,\ldots,0}_{i_0-2},-1,1,0,\ldots,0) &
                        if $1<i_0=j_0<\ell+1$.}
\]
We will often refer to this associated
element $m_p$ when we talk about a \dmove $p$.

Let us denote by $P_{ij}$ the set of moves from node $i$ on the left
to node $j$ on the right. For a move $p$, denote by $T_p$ the
corresponding  operator
on $\clh^{\otimes \ell(\omega)}\otimes \clk^{\otimes \ell}$.
Then
\be
\pi_\omega(u_{ij})=\sum_{p\in P_{ij}}T_p.
\ee
Denote by $W_p$ the operator obtained from
$T_p$ by replacing $\sqrt{I-q^{2N+2}}S$ by $S$
and $S^*\sqrt{I-q^{2N+2}}$ by $S^*$.
One can show easily that $m_p$ is the unique element
in $\bbz^\Lambda\times\bbz^{\Lambda_0}$ whose entries are all
$0$ or $\pm 1$ such that $\langle W_p e_\gamma,e_{\gamma+m_p}\rangle\neq 0$
for some $\gamma\in\Gamma$.

\blmma\label{movetemp}
Let $p, p'\in P_{ij}$. If $p$ and $p'$ are different, then
for some $(r,n)$, where $1\leq r\leq t$ and $n\in \Lambda_r$, one has
either $m_p(r,n)=0, m_{p'}(r,n)=\pm 1$
or $m_p(r,n)=\pm 1, m_{p'}(r,n)=0$.
\elmma
\prf
Since $p$ and $p'$ both belong to $P_{ij}$ and are different,
$m_p(r,n)\neq m_{p'}(r,n)$ for some pair $(r,n)$.
Now look at the coordinate where they are unequal for the first time
(from the left), that is, let $(r,n)$ be the pair such that
\[
r=\max \{1\leq j\leq t: m_p(j,i)\neq m_{p'}(j,i) \mbox{ for some }i\},\quad
n=\min \{i\in \Lambda_r: m_p(r,i)\neq m_{p'}(r,i)\}.
\]
It is easy to see now that for this pair $(r,n)$, the
required conclusion holds.\qed

\blmma
Let $F$ be a finite set of moves.
For $p\in F$, let $D_p$ be a (not necessarily bounded)
number operator, i.e.\ an operator
of the form $e_\gamma\mapsto t_\gamma e_\gamma$.
If $\sum_{p\in F}D_pW_p$ is bounded, then
$D_pW_p$ is bounded for each $p\in F$.
\elmma
\prf
Take $p'\in F$. Assume that $|F|>1$.
We will show that boundedness of $\sum_{p\in F}D_pW_p$
implies that of $\sum_{p\in F'}D_pW_p$
for some subset $F'$ of $F$ such that $p'\in F'$ and $|F'|<|F|$.

Let $p''\in F$ be an element of $F$ other than $p'$.
By the previous lemma, there is a pair $(r,n)$ such that
either $m_{p'}(r,n)=0$ and $m_{p''}(r,n)=\pm 1$
or $m_{p'}(r,n)=\pm 1$ and $m_{p''}(r,n)=0$.
For $z\in S^1$, let $U_z$ be the unitary operator on $L_2(\Gamma)$
given by $U_z e_\gamma=z^{\gamma(r,n)}e_\gamma$.
Now the proof will follow from the boundedness of the operator
$\int_{z\in S^1}U_z(\sum_{p\in F}D_pW_p)U_z^*\,dz$.\qed

\bppsn\label{bdd1}
$[D,\pi_\omega(u_{ij})]$ is bounded for all $i$ and $j$
if and only if
$[D,W_p]$ is bounded for all moves $p$.
\eppsn
\prf
It is enough to show that if $[D,\pi_\omega(u_{ij})]$
is bounded, and if $p\in P_{ij}$, then $[D,W_p]$ is bounded.
Since $\pi_\omega(u_{ij})=\sum_{p\in P_{ij}}T_p$
and each $[D, T_p]$ is of the form $D_pW_p$,
it follows from the forgoing lemma that each $[D, T_p]$
is bounded. Since $\sqrt{1-q^{2n+2}}$ is a bounded quantity
whose inverse is also bounded, it follows that
$[D, T_p]$ is bounded if and only if $[D,W_p]$ is bounded.\qed

Thus there is a positive constant $c$  such that
$D$ will have bounded commutators with all the $\pi_\omega(u_{ij})$'s
if and only if $\|[D,W_p]\|\leq c$.

Let $p=((i_t,j_t),\ldots,(i_0,j_0))$ be a move.
A coordinate $(r,s)$ is said to be a \textbf{diagonal component}
of $p$ if either $i_r<j_r$ and $s\in\{i_r,i_r+1,\ldots,j_r-1\}$,
or
$j_r=i_r-1=s$.
One can check that this would correspond exactly
to the diagonal parts of the move in the diagram representing $\omega$.
Denote by $c(\gamma,p)$ the quantity
$\sum_{(j,i)}\gamma(j,i)$, the sum being taken over all
diagonal components of $p$.

\blmma\label{bdd3}
$[D,W_p]$ is bounded if and only if
$|d(\gamma+m_p)-d(\gamma)|\leq cq^{-c(\gamma,p)}$.
\elmma
\prf
Follows easily once one writes down
the expression of the commutator.
\qed\\
An immediate corollary is the following.
\bcrlre
Let $H_i$ be as in definition~\ref{neqmove}. Then
$|d(\gamma+H_i)-d(\gamma)|\leq c$ for all $\gamma\in\Gamma$
and $1\leq i\leq \ell+1$.
\ecrlre

\newsection{The growth graph and sign characterization}
Let us now form the graph $\clg_c$ by connecting two vertices
$\gamma$ and $\gamma'$ if $|d(\gamma)-d(\gamma')|\leq c$.
Characterization of $\sgn D$ will then proceed as
outlined in the beginning of section~\ref{g_scheme}.

\bdfn
For $i\in \Lambda_0$, let $J_i$  be the set $\{j\geq 1:i\in \Lambda_j\}$.
The set
$\scrf =\{\gamma\in\bbz^\Lambda\times\bbz^{\Lambda_0}: -\gamma(0,i)=\gamma(0,i-1)=\gamma(j,i)
\mbox{ for all }j\in J_i\}$
will be called the \textbf{free plane}.
For a point $\gamma\in\Gamma$, we call the set
$\scrf_\gamma=\{\gamma+\gamma'\in\Gamma: \gamma'\in \scrf\}$
the \textbf{free plane passing through $\gamma$}.
\edfn
Note that for $\gamma\in\scrf$, the coordinates $\gamma(j,i)$
are all equal for $j\in J_i$.

For $1\leq i\leq \ell$, define
$j_i$ to be 0 if $J_i$ is empty, and to be that
element of $J_i$ for which
$\gamma(j_i,i)=\min \{\gamma(j,i):j\in J_i\}$.

\brmrk
1. If $J_i$ is nonempty, $j_i$ need not be unique.\\
2. If $\gamma'\in\scrf_\gamma$, then
$\min_j \gamma(j,i)$ and $\min_j \gamma'(j,i)$
are attained for the same set of values of $j$.
\ermrk

Note that given a $\gamma\in\Gamma$, elements in $\scrf_\gamma$ are determined
by the coordinates $(j_i,i)$, $i=1,\ldots,\ell$.

\blmma\label{newlemma}
Let $\gamma,\gamma'\in\Gamma$.
Then either $\scrf_\gamma=\scrf_{\gamma'}$
or $\scrf_\gamma$ and $\scrf_{\gamma'}$ are disjoint.
\elmma
\prf
Proof folows from the observation that
$\scrf_\gamma=\gamma+\scrf$ and $\scrf$ is a
subgroup in $\bbz^\Lambda\times\bbz^{\Lambda_0}$.
\qed

\blmma\label{sweep5}
Let $\gamma\in \Gamma$, and $\gamma'\in \scrf_\gamma$.
Let $\gamma''$ be the element in $\scrf_\gamma$ for which
\[
\gamma''(j_\ell,\ell)=\gamma'(j_\ell,\ell),\quad
  \gamma''(j_i,i)= 0 \mbox{ for all } i<\ell.
\]
Then there is a path in $\scrf_\gamma$ joining $\gamma'$ to $\gamma''$
such that throughout this path, the $(j_\ell,\ell)$-coordinate
remains constant.
\elmma
\prf
Apply successively the moves
\[
\gamma(j_{\ell-1},\ell-1)H_{\ell-1},\quad
 (\gamma(j_{\ell-2},\ell-2)+\gamma(j_{\ell-1},\ell-1))H_{\ell-2},\quad
\ldots,\quad
\left(\sum_{i=1}^{\ell-1} \gamma(j_i,i)\right) H_1.
\]
As none of these moves touch the $(j_\ell,\ell)$-coordinate,
it remains constant throughout the path.\qed

\blmma\label{signfree2}
Let $\gamma\in\Gamma$. Then either $\scrf_\gamma^+$ is finite
or $\scrf_\gamma^-$ is finite.
\elmma
\prf
Write $C(\gamma)=\gamma(j_\ell,\ell)$.
We will first show that $C(\scrf_\gamma^+)$
and $C(\scrf_\gamma^-)$ can not both be infinite.
Suppose if possible both  $C(\scrf_\gamma^+)$
and $C(\scrf_\gamma^-)$ are infinite. Then there
exist two sequences of elements
$\gamma_n$ and $\delta_n$ with
$\gamma_n\in\scrf_\gamma^+$
and
$\delta_n\in\scrf_\gamma^-$
such that
\[
C(\gamma_1)<C(\delta_1)<C(\gamma_2)<C(\delta_2)<
\cdots.
\]
Now start at $\gamma_n$ and employ lemma~\ref{sweep5}
to reach a point $\gamma_n'\in\scrf_\gamma^+$ such that
$C(\cdot)$ remains constant throughout the path
and for which
\[
\gamma_n'(j_i,i)=0 \mbox{ for all }i<\ell.
\]
Similarly start at $\delta_n$ and employ lemma~\ref{sweep5}
to use a path where $C(\cdot)$ remains constant
to reach a point $\delta_n'\in\scrf_\gamma^-$for which
\[
\delta_n'(j_i,i)=0 \mbox{ for all }i<\ell.
\]
Now use the move $H_\ell$ to go from $\gamma_n'$
to $\delta_n'$. The paths thus constructed
are all disjoint,
because the $C(\cdot)$ coordinate lies between
$C(\gamma_n)$ and $C(\delta_n)$ throughtout.
But that means $(\scrf_\gamma^+,\scrf_\gamma^-)$
admits an infinite ladder.

Next, suppose $C(\scrf_\gamma^-)\seq [-K,K]$.
If $\{\gamma'(j_i,i):\gamma'\in \scrf_\gamma\}$
is not bounded for some $i$ with $1\leq i\leq \ell-1$,
get a sequence of points $\gamma_n \in  \scrf_\gamma$
such that
$\gamma_n(j_i,i)<\gamma_{n+1}(j_i,i)$ for all $n$.
Starting at each $\gamma_n$, apply the move $H_{\ell+1}$
enough (e.g.\ $2K+1$) times to produce an infinite ladder.
\qed

Let us next define a set that will play the role of a complementary axis.
Let
\[
\scrc =\{\gamma\in\Gamma: \prod_{j\in J_i} \gamma(j,i)=0 \mbox{ for all }i\}.
\]
It follows from the sweepout argument used in the proof
of lemma~\ref{sweep5} that for any $\gamma'\in\Gamma$,
there is a $\gamma\in \scrc $
such that $\gamma'\in \scrf_\gamma$.
But it is not necessary that for two distinct elements
$\gamma$ and $\gamma'$ in $\scrc $, $\scrf_\gamma$ and $\scrf_{\gamma'}$ are disjoint.
However, this will not be of serious concern to us.

Let
\[
i_{min}=\min\{i\in \Lambda_0: |J_i|>1\},\quad
j_{min}=\min J_{i_{min}},\quad
j_{max}=\max J_{i_{min}}.
\]
Thus $i_{min}$ is the minimum $i$ for which $s_i$ appears more than once
in $\omega$, $j_{min}$ and $j_{max}$ are the first and the last string
where it appears.
Suppose now that we have removed the horizontal arrows labelled
$+$ or $-$ corresponding to all the $s_i$'s for which $|J_i|=1$.
Note that this would in particular remove all labelled horizontal
lines corresponding to $s_i$'s for $i<i_{min}$.
Suppose the $j_{min}$\raisebox{.4ex}{th} segment of a move is $(i_{min},i_{min})$.
This will uniquely specify the 0\raisebox{.4ex}{th} segment which will
be of the form $(i_0,i_0)$ for some $i_0\leq i_{min}$.
Now define
\be\label{defn_c0c1}
C_0(\gamma):=\gamma(j_{min},i_{min})+\gamma(0,i_0),\qquad
C_1(\gamma)=\gamma(j_{max},i_{min})
\ee
for $\gamma\in\Gamma$.

\blmma\label{sweep4}
Let $\gamma\in \scrc $. Define an element $\gamma'\in\Gamma$
by the following prescription:
\[
\gamma'(j,i)=0 \mbox{ for all } j\geq 1,\quad
\gamma'(0,i)=\cases{0 & if $i\neq i_0$,\cr
            C_0(\gamma) & if $i=i_0$.}
\]
Then there is a path connecting $\gamma$ to $\gamma'$ such
that $C_0(\cdot)$ remains constant throughout this path.
\elmma
\prf
We will describe a recursive algorithm to go
from $\gamma$ to $\gamma'$.
Observe that since $\gamma\in \scrc $, we have
$\gamma(t,\max \Lambda_t)=0$.
To begin with, remove all the horizontal arrows labelled
$+$ or $-$ corresponding to the $s_i$'s for which $|J_i|=1$,
and work with the resulting diagram.

Now suppose we are at $\delta\in\Gamma$ which satisfies
\[
\delta(j,i)=0 \mbox{ for all } j > r,\quad \delta(r,i)=0
  \mbox{ for all } i > n\in \Lambda_r.
\]
\textbf{Step I.}
\begin{quotation}
  \textbf{Case I}. $r=j_{min}$ and $n=i_{min}$: then apply the move whose
  $j_{min}$\raisebox{.4ex}{th} segment is $(i_{min},i_{min})$.  Apply
  this $\delta(j_{min},i_{min})$ times.  This will make the
  $(j_{min},i_{min})$-coordinate zero and the $(0,i_0)$-coordinate
  $C_0(\gamma)$.  Now proceed to step II.
\end{quotation}

\begin{quotation}
  \textbf{Case II}. $r\neq j_{min}$ or $n\neq i_{min}$:
Proceed with the following algorithm.\\[1ex]
\begin{tabular}{|p{400pt}|}
\hline
\textbf{Algorithm $\blda(r,n)$.} ($\min \Lambda_r\leq n \leq \max \Lambda_r$)
\par Remove all horizontal arrows labelled $+$ or $-$ from the $s_i$'s in
the strings $s_{\Lambda_t}$, $s_{\Lambda_{t-1}}$, $\ldots$, $s_{\Lambda_{r+1}}$ as well as from
the $s_i$'s corresponding to $i\in \Lambda_r$, $i> n$.
What this will achieve is the following:
any permissible move in the resulting diagram
will not change the coordinates
$(j,i)$ where either $r+1\leq j\leq t$ or $j=r$ and $i> n$.
\newline Apply the negative of the move  whose
$r$\raisebox{.4ex}{th} segment is
$(n+1,\max \Lambda_r+1)$
for $\delta(r,n)$ number of times.
This would kill the $(r,n)$-coordinate, i.e.\ will make it zero.
Now remove the two horizontal lines labelled `$+$' and `$-$'
corresponding to $s_n$ appearing in the string $s_{\Lambda_r}$.  \\
\hline
\end{tabular}
\end{quotation}

\vspace{1ex}

\noindent\textbf{Step II.}
\begin{quotation}
  \noindent\textbf{Case I.} $n>\min \Lambda_r$: keep $r$ intact, reduce the value
  of $n$ by 1 and go back to step I.\\
  \textbf{Case II.} $r>1$ and $n=\min \Lambda_r$: change $n$ to $\max
  \Lambda_{r-1}$, then
  reduce the value of $r$ by 1, and go back to step I.\\
  \textbf{Case III.}  $r=1$ and $n=\min \Lambda_1$: proceed to step III.
\end{quotation}
\textbf{Step III.}
All the $(j,i)$-coordinates for $j\geq 1$ are now zero.
Next, apply moves
ending at $i$ for
$i>i_0+1$ appropriate number of times starting from
the top to kill the coordinates $(0,i)$ for $i>i_0$.
Thus we have now reached an element $\delta$ for which
$\delta(j,i)=0$ whenever   $j\geq 1, i\in \Lambda_j$
    or $j=0, i>i_0$.
Therefore we now need to kill the coordinates
$(0,i)$ for $i<i_0$.
This is achieved as follows. Remove the
horizontal arrows labelled $+$ or $-$ from all $s_i$'s.
Now apply the moves
ending at $i$ for $i<i_0$ appropriate
number of times starting from the bottom.\qed

The next  diagram and the table that follows it
will explain the proof in a simple case.\\[2ex]

\hspace*{20pt}
\setlength{\unitlength}{0.00041667in}
\begingroup\makeatletter\ifx\SetFigFont\undefined%
\gdef\SetFigFont#1#2#3#4#5{%
  \reset@font\fontsize{#1}{#2pt}%
  \fontfamily{#3}\fontseries{#4}\fontshape{#5}%
  \selectfont}%
\fi\endgroup%
{\renewcommand{\dashlinestretch}{30}
\begin{picture}(11862,7512)(0,-10)
\path(1500,6306)(2700,6306)(3900,6306)
        (5100,6306)(6300,6306)(7500,6306)(8700,6306)
\path(1500,2706)(2700,2706)(3900,2706)
        (5100,2706)(6300,2706)(7500,2706)(8700,2706)
\path(1500,1506)(2700,1506)(3900,1506)
        (5100,1506)(6300,1506)(7500,1506)(8700,1506)
\path(1500,3906)(2700,3906)(3900,3906)
        (5100,3906)(6300,3906)(7500,3906)(8700,3906)
\path(1500,5106)(2700,5106)(3900,5106)
        (5100,5106)(6300,5106)(7500,5106)(8700,5106)
\path(1500,2706)(2700,3906)(3900,5106)(5100,6306)
\path(1500,3906)(2700,2706)
\path(2700,5106)(3900,3906)
\path(3900,6306)(5100,5106)
\path(5100,5106)(6300,3906)(7500,2706)(8700,1506)
\path(5100,3906)(6300,5106)
\path(6300,2706)(7500,3906)
\path(7500,1506)(8700,2706)
\dashline{60.000}(1500,4056)(2700,4056)(5100,6456)(8700,6456)
\dashline{60.000}(1500,5256)(3675,5256)(5100,6606)(8700,6606)
\dashline{60.000}(1500,3756)(2700,2556)(7500,2556)(8700,1356)
\dashline{60.000}(1500,6456)(3900,6456)(5100,5256)(8625,5256)
\path(900,6981)(1725,6531)
\path(1605.287,6562.125)(1725.000,6531.000)(1634.018,6614.799)
\path(900,5781)(1725,5331)
\path(1605.287,5362.125)(1725.000,5331.000)(1634.018,5414.799)
\path(900,4581)(1725,4131)
\path(1605.287,4162.125)(1725.000,4131.000)(1634.018,4214.799)
\path(975,3081)(1650,3531)
\path(1566.795,3439.474)(1650.000,3531.000)(1533.513,3489.397)
\path(8775,1506)(9450,1506)
\path(9330.000,1476.000)(9450.000,1506.000)(9330.000,1536.000)
\path(9600,1506)(10275,1506)
\path(10155.000,1476.000)(10275.000,1506.000)(10155.000,1536.000)
\path(11175,1506)(11850,1506)
\path(11730.000,1476.000)(11850.000,1506.000)(11730.000,1536.000)
\path(10350,1506)(11025,1506)
\path(10905.000,1476.000)(11025.000,1506.000)(10905.000,1536.000)
\path(8775,2706)(9450,2706)
\path(9330.000,2676.000)(9450.000,2706.000)(9330.000,2736.000)
\path(9600,2706)(10275,2706)
\path(10155.000,2676.000)(10275.000,2706.000)(10155.000,2736.000)
\path(11175,2706)(11850,2706)
\path(11730.000,2676.000)(11850.000,2706.000)(11730.000,2736.000)
\path(10350,2706)(11025,2706)
\path(10905.000,2676.000)(11025.000,2706.000)(10905.000,2736.000)
\path(8775,3906)(9450,3906)
\path(9330.000,3876.000)(9450.000,3906.000)(9330.000,3936.000)
\path(9600,3906)(10275,3906)
\path(10155.000,3876.000)(10275.000,3906.000)(10155.000,3936.000)
\path(11175,3906)(11850,3906)
\path(11730.000,3876.000)(11850.000,3906.000)(11730.000,3936.000)
\path(10350,3906)(11025,3906)
\path(10905.000,3876.000)(11025.000,3906.000)(10905.000,3936.000)
\path(8775,5106)(9450,5106)
\path(9330.000,5076.000)(9450.000,5106.000)(9330.000,5136.000)
\path(9600,5106)(10275,5106)
\path(10155.000,5076.000)(10275.000,5106.000)(10155.000,5136.000)
\path(11175,5106)(11850,5106)
\path(11730.000,5076.000)(11850.000,5106.000)(11730.000,5136.000)
\path(10350,5106)(11025,5106)
\path(10905.000,5076.000)(11025.000,5106.000)(10905.000,5136.000)
\path(8775,6306)(9450,6306)
\path(9330.000,6276.000)(9450.000,6306.000)(9330.000,6336.000)
\path(9600,6306)(10275,6306)
\path(10155.000,6276.000)(10275.000,6306.000)(10155.000,6336.000)
\path(11175,6306)(11850,6306)
\path(11730.000,6276.000)(11850.000,6306.000)(11730.000,6336.000)
\path(10350,6306)(11025,6306)
\path(10905.000,6276.000)(11025.000,6306.000)(10905.000,6336.000)
\dashline{60.000}(8700,1356)(11850,1356)
\dashline{60.000}(8625,5256)(11850,5256)
\dashline{60.000}(8700,6456)(11850,6456)
\dashline{60.000}(8700,6606)(11850,6606)
\path(1620.000,7386.000)(1500.000,7356.000)(1620.000,7326.000)
\path(1500,7356)(4500,7356)
\path(5625,7356)(8700,7356)
\path(8580.000,7326.000)(8700.000,7356.000)(8580.000,7386.000)
\path(8970.000,7386.000)(8850.000,7356.000)(8970.000,7326.000)
\path(8850,7356)(9600,7356)
\path(10650,7356)(11850,7356)
\path(11730.000,7326.000)(11850.000,7356.000)(11730.000,7386.000)
\dashline{60.000}(1500,2556)(5175,6231)(11850,6231)
\dashline{60.000}(1500,4956)(2700,4956)(3900,3756)
        (4650,3756)(4725,3756)(6450,5481)(11775,5481)
\dashline{60.000}(1500,3531)(2700,2331)(6150,2331)
        (7575,3756)(11850,3756)
\path(611,5427)(1436,4977)
\path(1316.287,5008.125)(1436.000,4977.000)(1345.018,5060.799)
\path(825,2106)(1500,2556)
\path(1416.795,2464.474)(1500.000,2556.000)(1383.513,2514.397)
\path(609,3997)(1434,3547)
\path(1314.287,3578.125)(1434.000,3547.000)(1343.018,3630.799)
\put(300,6906){\makebox(0,0)[lb]{\smash{{{\SetFigFont{9}{10.8}{\rmdefault}{\mddefault}{\updefault}$m_3$}}}}}
\put(300,4506){\makebox(0,0)[lb]{\smash{{{\SetFigFont{9}{10.8}{\rmdefault}{\mddefault}{\updefault}$m_2$}}}}}
\put(300,3006){\makebox(0,0)[lb]{\smash{{{\SetFigFont{9}{10.8}{\rmdefault}{\mddefault}{\updefault}$m_4$}}}}}
\put(300,5706){\makebox(0,0)[lb]{\smash{{{\SetFigFont{9}{10.8}{\rmdefault}{\mddefault}{\updefault}$m_1$}}}}}
\put(2175,81){\makebox(0,0)[lb]{\smash{{{\SetFigFont{9}{10.8}{\rmdefault}{\mddefault}{\updefault}$t=4$, $i_{min}=2$, $j_{min}=2$, $j_{max}=4$, $i_0=1$}}}}}
\put(4350,6381){\makebox(0,0)[lb]{\smash{{{\SetFigFont{9}{10.8}{\rmdefault}{\mddefault}{\updefault}$+$}}}}}
\put(3075,4881){\makebox(0,0)[lb]{\smash{{{\SetFigFont{9}{10.8}{\rmdefault}{\mddefault}{\updefault}$+$}}}}}
\put(1875,3681){\makebox(0,0)[lb]{\smash{{{\SetFigFont{9}{10.8}{\rmdefault}{\mddefault}{\updefault}$+$}}}}}
\put(5475,4881){\makebox(0,0)[lb]{\smash{{{\SetFigFont{9}{10.8}{\rmdefault}{\mddefault}{\updefault}$+$}}}}}
\put(6675,3981){\makebox(0,0)[lb]{\smash{{{\SetFigFont{9}{10.8}{\rmdefault}{\mddefault}{\updefault}$+$}}}}}
\put(7875,2781){\makebox(0,0)[lb]{\smash{{{\SetFigFont{9}{10.8}{\rmdefault}{\mddefault}{\updefault}$+$}}}}}
\put(4275,4881){\makebox(0,0)[lb]{\smash{{{\SetFigFont{9}{10.8}{\rmdefault}{\mddefault}{\updefault}$-$}}}}}
\put(3075,3681){\makebox(0,0)[lb]{\smash{{{\SetFigFont{9}{10.8}{\rmdefault}{\mddefault}{\updefault}$-$}}}}}
\put(1875,2481){\makebox(0,0)[lb]{\smash{{{\SetFigFont{9}{10.8}{\rmdefault}{\mddefault}{\updefault}$-$}}}}}
\put(7950,1281){\makebox(0,0)[lb]{\smash{{{\SetFigFont{9}{10.8}{\rmdefault}{\mddefault}{\updefault}$-$}}}}}
\put(6750,2781){\makebox(0,0)[lb]{\smash{{{\SetFigFont{9}{10.8}{\rmdefault}{\mddefault}{\updefault}$-$}}}}}
\put(5475,3681){\makebox(0,0)[lb]{\smash{{{\SetFigFont{9}{10.8}{\rmdefault}{\mddefault}{\updefault}$-$}}}}}
\put(8850,1581){\makebox(0,0)[lb]{\smash{{{\SetFigFont{9}{10.8}{\rmdefault}{\mddefault}{\updefault}$+$}}}}}
\put(9675,2781){\makebox(0,0)[lb]{\smash{{{\SetFigFont{9}{10.8}{\rmdefault}{\mddefault}{\updefault}$+$}}}}}
\put(10425,3981){\makebox(0,0)[lb]{\smash{{{\SetFigFont{9}{10.8}{\rmdefault}{\mddefault}{\updefault}$+$}}}}}
\put(8925,2781){\makebox(0,0)[lb]{\smash{{{\SetFigFont{9}{10.8}{\rmdefault}{\mddefault}{\updefault}$-$}}}}}
\put(9750,3981){\makebox(0,0)[lb]{\smash{{{\SetFigFont{9}{10.8}{\rmdefault}{\mddefault}{\updefault}$-$}}}}}
\put(11325,6006){\makebox(0,0)[lb]{\smash{{{\SetFigFont{9}{10.8}{\rmdefault}{\mddefault}{\updefault}$-$}}}}}
\put(11250,4806){\makebox(0,0)[lb]{\smash{{{\SetFigFont{9}{10.8}{\rmdefault}{\mddefault}{\updefault}$+$}}}}}
\put(10500,4806){\makebox(0,0)[lb]{\smash{{{\SetFigFont{9}{10.8}{\rmdefault}{\mddefault}{\updefault}$-$}}}}}
\put(4650,7281){\makebox(0,0)[lb]{\smash{{{\SetFigFont{9}{10.8}{\rmdefault}{\mddefault}{\updefault}$\mathbb{N}$ part}}}}}
\put(9675,7281){\makebox(0,0)[lb]{\smash{{{\SetFigFont{9}{10.8}{\rmdefault}{\mddefault}{\updefault}$\mathbb{Z}$ part}}}}}
\put(2175,681){\makebox(0,0)[lb]{\smash{{{\SetFigFont{9}{10.8}{\rmdefault}{\mddefault}{\updefault}$\omega=(s_2s_3s_4)(s_3)(s_2)(s_1)$}}}}}
\put(225,2031){\makebox(0,0)[lb]{\smash{{{\SetFigFont{9}{10.8}{\rmdefault}{\mddefault}{\updefault}$m_7$}}}}}
\put(0,3906){\makebox(0,0)[lb]{\smash{{{\SetFigFont{9}{10.8}{\rmdefault}{\mddefault}{\updefault}$m_5$}}}}}
\put(0,5331){\makebox(0,0)[lb]{\smash{{{\SetFigFont{9}{10.8}{\rmdefault}{\mddefault}{\updefault}$m_6$}}}}}
\end{picture}
}

\vspace{4ex}

The table below illustrates the sweepout procedure
described in the proof of lemma~\ref{sweep4}.
Starting from a point $\gamma\in\scrc$, it shows the successive moves
applied and  how the resulting element looks like
at each stage. Observe that for any $\gamma\in\scrc$, one must have
$\gamma(4,4)=0=\gamma(1,1)$.\\[2ex]
\hspace*{-2mm}
\begin{tabular}{|r||r|r|r|r|r|r||r|r|r|r|}
\hline
coordinate  & (4,2) & (4,3) & (4,4) & (3,3) & (2,2) & (1,1) & (0,1) & (0,2) & (0,3) & (0,4) \\
\hline
\hline
$\gamma$ & $\ast$ &  $\ast$ &  0 & $\ast$  & $a$ & 0 & $b$ & $\ast$  & $\ast$  & $\ast$ \\
\hline
move $m_1$ & 0 & $+1$ & 0 & 0 & 0 & 0 & 0 & 0 & 0 & $-1$ \\
\hline
$\gamma_1=-\gamma(4,3)m_1(\gamma)$
          & $\ast$ &  0 &  0 & $\ast$  & $a$ & 0 & $b$ & $\ast$  & $\ast$  & $\ast$ \\
\hline
\hline
move $m_2$ & $+1$ & 0 & 0 & 0 & 0 & 0 & 0 & 0 & 0 & $-1$ \\
\hline
$\gamma_2=-\gamma(4,2)m_2(\gamma_1)$
          & 0 &  0 &  0 & $\ast$  & $a$ & 0 & $b$ & $\ast$  & $\ast$  & $\ast$ \\
\hline
\hline
move $m_3$ & 0 & 0 & 0 & $+1$ & 0 & 0 & 0 & 0 & $-1$ & $+1$ \\
\hline
$\gamma_3=-\gamma(3,3)m_3(\gamma_2)$
          & 0 & 0 &  0 & 0  & $a$ & 0 & $b$ & $\ast$  & $\ast$  & $\ast$ \\
\hline
\hline
move $m_4$ & 0 & 0 & 0 & 0 & $-1$ & 0 & $+1$ & 0 & 0 &0  \\
\hline
$\gamma_4=\gamma(2,2)m_4(\gamma_3)$
      & 0 &  0 & 0 & 0  & 0 & 0 & $a+b$ & $\ast$ & $\ast$  & $\ast$   \\
\hline
\hline
move $m_5$ & 0 & 0 & 0 & 0 & 0 & 0 & 0 & $-1$ & $+1$ &0  \\
\hline
$\gamma_5=\gamma_4(0,2)m_5(\gamma_4)$
      & 0 &  0 & 0 & 0  & 0 & 0 & $a+b$ & 0 & $\ast$  & $\ast$   \\
\hline
\hline
move $m_6$ & 0 & 0 & 0 & 0 & 0  & 0 & 0 & 0 & $-1$ &$+1$  \\
\hline
$\gamma_6=\gamma_5(0,3)m_6(\gamma_5)$
      & 0 &  0 & 0 & 0  & 0 & 0 & $a+b$ & 0 & 0 & $\ast$   \\
\hline
\hline
move $m_7$ & 0 & 0 & 0 & 0 & 0 & 0 & 0 & 0 & 0 & $-1$  \\
\hline
$\gamma'=\gamma_6(0,4)m_7(\gamma_6)$
      & 0 &  0 & 0 & 0  & 0 & 0 & $a+b$ & 0 & 0  & 0   \\
\hline
\end{tabular}

\blmma\label{about_c0}
Both $C_0(\scrc ^+)$ and $C_0(\scrc ^-)$ can not be infinite.
\elmma
\prf
If both are infinite, there would exist elements
$\gamma_n\in \scrc ^+$ and $\delta_n\in \scrc ^-$ such that
\[
C_0(\gamma_1)<C_0(\delta_1)<C_0(\gamma_2)<C_0(\delta_2)<\ldots.
\]
Let $\gamma_n'$ and $\delta_n'$ be given by
\[
\gamma_n'(j,i)=0 \mbox{ for all } j\geq 1,\quad
\gamma_n'(0,i)=\cases{0 & if $i\neq i_0$,\cr
         C_0(\gamma_n) & if $i=i_0$,}
\]
\[
\delta_n'(j,i)=0 \mbox{ for all } j\geq 1,\quad
\delta_n'(0,i)=\cases{0 & if $i\neq i_0$,\cr
             C_0(\delta_n) & if $i=i_0$.}
\]
Use the earlier lemma to get paths between $\gamma_n$ and $\gamma_n'$
and between $\delta_n$ and $\delta_n'$.
Remove all the labelled arrows from  all the $s_i$'s.
Let $m_i$ be the move in the resulting diagram whose
0\raisebox{.4ex}{th} segment is $(i,i)$,
and let $m=\sum_{{i=1}\atop\leftarrow}^{i_0}m_i$.
Apply this move $C_0(\delta_n)-C_0(\gamma_n)$ times
to connect $\gamma_n'$ and $\delta_n'$.
Thus there is a path $p_n$ connecting $\gamma_n$ and $\delta_n$,
and throughout this path, $C_0(\cdot)$ lies between $C_0(\gamma_n)$
and $C_0(\delta_n)$. Therefore the paths $p_n$ are disjoint.\qed

For the next two lemmas, we will assume
that  $C_0(\scrc ^-)$ is finite.

\blmma
Assume  $C_0(\scrc ^-)$ is finite.
Let $C_1$ be as defined prior to lemma~\ref{sweep4},
i.e.\ $C_1(\gamma)=\gamma(j_{max},i_{min})$.
Then the set $C_1(\scrc ^-)$ is finite.
\elmma
\prf
Let $K\in\bbn$ be such that $C_0(\scrc ^-)\seq [-K,K]$.
If  $C_1(\scrc ^-)$ is not finite, there is a $\gamma_n\in \scrc ^-$ such that
\[
C_1(\gamma_1) < C_1(\gamma_2) <C_1(\gamma_3) < \ldots.
\]
Now the idea is to get a path $p_n$ joining $\gamma_n$ to some
$\delta_n$ such that $C_1(\cdot)$ remains constant throughout $p_n$,
and $C_0(\delta_n)>K$, so that each $\delta_n\in \scrc ^+$.

Start at $\gamma_n$. Apply algorithm $\blda(r,n)$
for
\[
\brray{lcll}
r&=&t,t-1,\ldots,j_{max}+1,&  \min \Lambda_r\leq n\leq \max \Lambda_r,\\
    r&=&j_{max},& i_{min}+1\leq n,\\
    r &<& j_{max},& \min \Lambda_r\leq n\leq \max \Lambda_r.
\erray
\]
Now apply the move $m_{i_0}$, where $m_i$'s are the moves
described in the proof of the previous
lemma, $3K$ times.\qed

Again we give a diagram  and a table to illustrate  the above proof
for the case $\omega=(s_2s_3s_4)(s_3)(s_2)(s_1)$.\\[2ex]

\hspace*{20pt}
\setlength{\unitlength}{0.00041667in}
\begingroup\makeatletter\ifx\SetFigFont\undefined%
\gdef\SetFigFont#1#2#3#4#5{%
  \reset@font\fontsize{#1}{#2pt}%
  \fontfamily{#3}\fontseries{#4}\fontshape{#5}%
  \selectfont}%
\fi\endgroup%
{\renewcommand{\dashlinestretch}{30}
\begin{picture}(11787,7512)(0,-10)
\path(1425,6306)(2625,6306)(3825,6306)
        (5025,6306)(6225,6306)(7425,6306)(8625,6306)
\path(1425,2706)(2625,2706)(3825,2706)
        (5025,2706)(6225,2706)(7425,2706)(8625,2706)
\path(1425,1506)(2625,1506)(3825,1506)
        (5025,1506)(6225,1506)(7425,1506)(8625,1506)
\path(1425,3906)(2625,3906)(3825,3906)
        (5025,3906)(6225,3906)(7425,3906)(8625,3906)
\path(1425,5106)(2625,5106)(3825,5106)
        (5025,5106)(6225,5106)(7425,5106)(8625,5106)
\path(1425,2706)(2625,3906)(3825,5106)(5025,6306)
\path(1425,3906)(2625,2706)
\path(2625,5106)(3825,3906)
\path(3825,6306)(5025,5106)
\path(5025,5106)(6225,3906)(7425,2706)(8625,1506)
\path(5025,3906)(6225,5106)
\path(6225,2706)(7425,3906)
\path(7425,1506)(8625,2706)
\dashline{60.000}(1425,5256)(3600,5256)(5025,6606)(8625,6606)
\dashline{60.000}(1425,6456)(3825,6456)(5025,5256)(8550,5256)
\path(8700,1506)(9375,1506)
\path(9255.000,1476.000)(9375.000,1506.000)(9255.000,1536.000)
\path(9525,1506)(10200,1506)
\path(10080.000,1476.000)(10200.000,1506.000)(10080.000,1536.000)
\path(11100,1506)(11775,1506)
\path(11655.000,1476.000)(11775.000,1506.000)(11655.000,1536.000)
\path(10275,1506)(10950,1506)
\path(10830.000,1476.000)(10950.000,1506.000)(10830.000,1536.000)
\path(8700,2706)(9375,2706)
\path(9255.000,2676.000)(9375.000,2706.000)(9255.000,2736.000)
\path(9525,2706)(10200,2706)
\path(10080.000,2676.000)(10200.000,2706.000)(10080.000,2736.000)
\path(11100,2706)(11775,2706)
\path(11655.000,2676.000)(11775.000,2706.000)(11655.000,2736.000)
\path(10275,2706)(10950,2706)
\path(10830.000,2676.000)(10950.000,2706.000)(10830.000,2736.000)
\path(8700,3906)(9375,3906)
\path(9255.000,3876.000)(9375.000,3906.000)(9255.000,3936.000)
\path(9525,3906)(10200,3906)
\path(10080.000,3876.000)(10200.000,3906.000)(10080.000,3936.000)
\path(11100,3906)(11775,3906)
\path(11655.000,3876.000)(11775.000,3906.000)(11655.000,3936.000)
\path(10275,3906)(10950,3906)
\path(10830.000,3876.000)(10950.000,3906.000)(10830.000,3936.000)
\path(8700,5106)(9375,5106)
\path(9255.000,5076.000)(9375.000,5106.000)(9255.000,5136.000)
\path(9525,5106)(10200,5106)
\path(10080.000,5076.000)(10200.000,5106.000)(10080.000,5136.000)
\path(11100,5106)(11775,5106)
\path(11655.000,5076.000)(11775.000,5106.000)(11655.000,5136.000)
\path(10275,5106)(10950,5106)
\path(10830.000,5076.000)(10950.000,5106.000)(10830.000,5136.000)
\path(8700,6306)(9375,6306)
\path(9255.000,6276.000)(9375.000,6306.000)(9255.000,6336.000)
\path(9525,6306)(10200,6306)
\path(10080.000,6276.000)(10200.000,6306.000)(10080.000,6336.000)
\path(11100,6306)(11775,6306)
\path(11655.000,6276.000)(11775.000,6306.000)(11655.000,6336.000)
\path(10275,6306)(10950,6306)
\path(10830.000,6276.000)(10950.000,6306.000)(10830.000,6336.000)
\dashline{60.000}(8550,5256)(11775,5256)
\dashline{60.000}(8625,6606)(11775,6606)
\path(1545.000,7386.000)(1425.000,7356.000)(1545.000,7326.000)
\path(1425,7356)(4425,7356)
\path(5550,7356)(8625,7356)
\path(8505.000,7326.000)(8625.000,7356.000)(8505.000,7386.000)
\path(8895.000,7386.000)(8775.000,7356.000)(8895.000,7326.000)
\path(8775,7356)(9525,7356)
\path(10575,7356)(11775,7356)
\path(11655.000,7326.000)(11775.000,7356.000)(11655.000,7386.000)
\dashline{60.000}(1425,6606)(3825,6606)(6675,3756)(11775,3756)
\dashline{60.000}(1425,6156)(3825,6156)(8625,1356)(11775,1356)
\path(675,4881)(1425,5256)
\path(1331.085,5175.502)(1425.000,5256.000)(1304.252,5229.167)
\path(660,5736)(1410,6111)
\path(1316.085,6030.502)(1410.000,6111.000)(1289.252,6084.167)
\path(585,6036)(1335,6411)
\path(1241.085,6330.502)(1335.000,6411.000)(1214.252,6384.167)
\path(600,6831)(1350,6606)
\path(1226.440,6611.747)(1350.000,6606.000)(1243.681,6669.217)
\put(2100,81){\makebox(0,0)[lb]{\smash{{{\SetFigFont{9}{10.8}{\rmdefault}{\mddefault}{\updefault}$t=4$, $i_{min}=2$, $j_{min}=2$, $j_{max}=4$, $i_0=1$}}}}}
\put(4275,6381){\makebox(0,0)[lb]{\smash{{{\SetFigFont{9}{10.8}{\rmdefault}{\mddefault}{\updefault}$+$}}}}}
\put(3000,4881){\makebox(0,0)[lb]{\smash{{{\SetFigFont{9}{10.8}{\rmdefault}{\mddefault}{\updefault}$+$}}}}}
\put(1800,3681){\makebox(0,0)[lb]{\smash{{{\SetFigFont{9}{10.8}{\rmdefault}{\mddefault}{\updefault}$+$}}}}}
\put(5400,4881){\makebox(0,0)[lb]{\smash{{{\SetFigFont{9}{10.8}{\rmdefault}{\mddefault}{\updefault}$+$}}}}}
\put(6600,3981){\makebox(0,0)[lb]{\smash{{{\SetFigFont{9}{10.8}{\rmdefault}{\mddefault}{\updefault}$+$}}}}}
\put(7800,2781){\makebox(0,0)[lb]{\smash{{{\SetFigFont{9}{10.8}{\rmdefault}{\mddefault}{\updefault}$+$}}}}}
\put(4200,4881){\makebox(0,0)[lb]{\smash{{{\SetFigFont{9}{10.8}{\rmdefault}{\mddefault}{\updefault}$-$}}}}}
\put(3000,3681){\makebox(0,0)[lb]{\smash{{{\SetFigFont{9}{10.8}{\rmdefault}{\mddefault}{\updefault}$-$}}}}}
\put(1800,2481){\makebox(0,0)[lb]{\smash{{{\SetFigFont{9}{10.8}{\rmdefault}{\mddefault}{\updefault}$-$}}}}}
\put(7875,1281){\makebox(0,0)[lb]{\smash{{{\SetFigFont{9}{10.8}{\rmdefault}{\mddefault}{\updefault}$-$}}}}}
\put(6675,2781){\makebox(0,0)[lb]{\smash{{{\SetFigFont{9}{10.8}{\rmdefault}{\mddefault}{\updefault}$-$}}}}}
\put(5400,3681){\makebox(0,0)[lb]{\smash{{{\SetFigFont{9}{10.8}{\rmdefault}{\mddefault}{\updefault}$-$}}}}}
\put(8775,1581){\makebox(0,0)[lb]{\smash{{{\SetFigFont{9}{10.8}{\rmdefault}{\mddefault}{\updefault}$+$}}}}}
\put(9600,2781){\makebox(0,0)[lb]{\smash{{{\SetFigFont{9}{10.8}{\rmdefault}{\mddefault}{\updefault}$+$}}}}}
\put(10350,3981){\makebox(0,0)[lb]{\smash{{{\SetFigFont{9}{10.8}{\rmdefault}{\mddefault}{\updefault}$+$}}}}}
\put(8850,2781){\makebox(0,0)[lb]{\smash{{{\SetFigFont{9}{10.8}{\rmdefault}{\mddefault}{\updefault}$-$}}}}}
\put(9675,3981){\makebox(0,0)[lb]{\smash{{{\SetFigFont{9}{10.8}{\rmdefault}{\mddefault}{\updefault}$-$}}}}}
\put(11250,6006){\makebox(0,0)[lb]{\smash{{{\SetFigFont{9}{10.8}{\rmdefault}{\mddefault}{\updefault}$-$}}}}}
\put(11175,4806){\makebox(0,0)[lb]{\smash{{{\SetFigFont{9}{10.8}{\rmdefault}{\mddefault}{\updefault}$+$}}}}}
\put(10425,4806){\makebox(0,0)[lb]{\smash{{{\SetFigFont{9}{10.8}{\rmdefault}{\mddefault}{\updefault}$-$}}}}}
\put(4575,7281){\makebox(0,0)[lb]{\smash{{{\SetFigFont{9}{10.8}{\rmdefault}{\mddefault}{\updefault}$\mathbb{N}$ part}}}}}
\put(9600,7281){\makebox(0,0)[lb]{\smash{{{\SetFigFont{9}{10.8}{\rmdefault}{\mddefault}{\updefault}$\mathbb{Z}$ part}}}}}
\put(2100,681){\makebox(0,0)[lb]{\smash{{{\SetFigFont{9}{10.8}{\rmdefault}{\mddefault}{\updefault}$\omega=(s_2s_3s_4)(s_3)(s_2)(s_1)$}}}}}
\put(75,4731){\makebox(0,0)[lb]{\smash{{{\SetFigFont{8}{9.6}{\rmdefault}{\mddefault}{\updefault}$m_1$}}}}}
\put(75,5556){\makebox(0,0)[lb]{\smash{{{\SetFigFont{8}{9.6}{\rmdefault}{\mddefault}{\updefault}$m_4$}}}}}
\put(0,5931){\makebox(0,0)[lb]{\smash{{{\SetFigFont{8}{9.6}{\rmdefault}{\mddefault}{\updefault}$m_2$}}}}}
\put(0,6756){\makebox(0,0)[lb]{\smash{{{\SetFigFont{8}{9.6}{\rmdefault}{\mddefault}{\updefault}$m_3$}}}}}
\end{picture}
}

\vspace{2ex}

The next table illustrates the argument in the above proof.
Starting from a point $\gamma\in\scrc^-$, it shows the successive moves
applied and  how the resulting element looks like
at each stage.\\[2ex]
\hspace*{-10mm}
\begin{tabular}{|r||r|r|r|r|r|r||r|r|r|r|}
\hline
coordinate  & (4,2) & (4,3) & (4,4) & (3,3) & (2,2) & (1,1) & (0,1) & (0,2) & (0,3) & (0,4) \\
\hline
\hline
$\gamma$ & $a$ &  $\ast$ &  0 & $\ast$  & $\ast$ & 0 & $b$ & $\ast$  & $\ast$  & $\ast$ \\
\hline
move $m_1$ & 0 & $+1$ & 0 & 0 & 0 & 0 & 0 & 0 & 0 & $-1$ \\
\hline
$\gamma_1=-\gamma(4,3)m_1(\gamma)$
          & $a$ &  0 &  0 & $\ast$  & $\ast$ & 0 & $b$ & $\ast$  & $\ast$  & $\ast$ \\
\hline
\hline
move $m_2$ & 0 & 0 & 0 & $+1$ & 0 & 0 & 0 & 0 &  $-1$ & $+1$\\
\hline
$\gamma_2=-\gamma(3,3)m_2(\gamma_1)$
          & $a$ &  0 &  0 & 0  & $\ast$ & 0 & $b$ & $\ast$  & $\ast$  & $\ast$ \\
\hline
\hline
move $m_3$ & 0 & 0 & 0 & 0& $+1$ & 0 & 0 &   $-1$ & $+1$&0 \\
\hline
$\gamma_3=-\gamma(2,2)m_3(\gamma_2)$
  & $a$ & 0 &  0 & 0  & 0 & 0 & $b$ & $\ast$  & $\ast$  & $\ast$ \\
\hline
\hline
move $m_4$ & 0 & 0 & 0 & 0 & 0 & 0 & $+1$ & 0 & 0 &0  \\
\hline
$\gamma_4=3Km_4(\gamma_3)$
      & $a$ &  0 & 0 & 0  & 0 & 0 & $b+3K$ & $\ast$ & $\ast$  & $\ast$   \\
\hline
\end{tabular}

\vspace{2ex}

\blmma
Assume $C_0(\scrc ^-)$ is finite.
Let
$C\equiv(j,i)$ be any coordinate other than $C_1\equiv (j_{max},i_{min})$.
Then $C(\scrc ^-)$ is finite.
\elmma
\prf
By the previous lemma, $C_1(\scrc ^-)$
is also bounded.
Let $K\in\bbn$ be such that $C_0(\scrc ^-)\seq [-K,K]$
and $C_1(\scrc ^-)\seq [-K,K]$.
The strategy would be the same as in the proof of the earlier
lemma with a slight modification.
If $C(\scrc ^-)$ is infinite,  we can choose $\gamma_n\in \scrc ^-$ such that
\[
C(\gamma_n)+K+1 < C(\gamma_{n+1})
\]
for every $n\in\bbn$.
Now connect every $\gamma_n$ to an element $\delta_n\in \scrc ^+$
by a path $p_n$ such that on $p_n$, the $C_1$ coordinate does not vary by
more than $K$. This will ensure that the paths $p_n$ are all disjoint.

For getting $p_n$ as described above, start at $\gamma_n$ and apply
successively the moves
\[
H_{\ell+1}, H_\ell, \ldots, H_{i_{min}+1},
\]
each one $K+1$ times. This will increase the $C_1$-coordinate by
$K+1$. Therefore the endpoint of the path will lie in $\scrc^+$.
\qed

Thus it now follows that
if $C_0(\scrc^-)$ is finite, then $\scrc ^-$ is finite.
Similar argument would  tell us that
if $C_0(\scrc^+)$ is finite, then $\scrc ^+$ is finite.
Therefore by lemma~\ref{about_c0}, either $\scrc^+$ or $\scrc^-$
is finite.
This, together with
proposition~\ref{signfree2} will give us the following theorem.
\bthm\label{noneqsign}
Let $D$ be a Dirac operator on $L_2(\Gamma)$
that diagonalises with respect to the canonical orthonormal basis.
Then up to a compact perturbation,
$\sgn D$ must be of the form $2P-I$ or $I-2P$
where $P$ is a projection onto the closed linear span of
$\{e_\gamma: \gamma\in\cup_{i=1}^k \scrf_{\gamma_i}\}$ for some
finite collection $\gamma_1,\gamma_2,\ldots,\gamma_k$ in $\Gamma$.
\ethm
\prf
We  first claim that there are finitely many
free planes $\scrf_\gamma$ for which
both $\scrf_\gamma^+$ and $\scrf_\gamma^-$
are nonempty.
It follows from the proofs of lemma~\ref{sweep5}
and proposition~\ref{signfree2} that any two points
on a free plane can be connected by a path
lying entirely on that plane. Therefore,
if there are infinitely many distinct $\scrf_\gamma$'s
for which  both $\scrf_\gamma^+$ and $\scrf_\gamma^-$
are nonempty, one can easily produce an infinite ladder.
Thus the claim is established.

Since for each $\gamma$, either $\scrf_\gamma^+$
or $\scrf_\gamma^-$ is finite, employing an appropriate
compact perturbation, it is now possible to ensure
that for every $\gamma\in\Gamma$, either
$\scrf_\gamma\seq \Gamma^+$
or
$\scrf_\gamma\seq \Gamma^-$.
This, along with the fact that either $\scrc^+$
or $\scrc^-$ must be finite, imply the desired result.\qed

We next show that under this restriction, compactness of the
commutator $[\sgn\, D, u_{ij}]$, or, equivalently, that of
$[P,u_{ij}]$'s will imply that
 $\sgn\, D$ is trivial.

Let $\gamma_1,\gamma_2,\ldots,\gamma_k$ be elements in $\Gamma$
and let $P$ be the projection onto
$\mbox{span}\,\{e_\gamma: \gamma\in \cup_i\scrf_{\gamma_i}\}$.
Then for any operator $T$,
we have
\[
[P,T]e_\gamma=\cases{PT e_\gamma & if $\gamma\not\in \cup_i\scrf_{\gamma_i}$,\cr
   (P-I)Te_\gamma & if $\gamma\in \cup_i\scrf_{\gamma_i}$.}
\]
Now let $r=\max \Lambda_t$ and take $T=\pi_\omega(u_{r+1,r})$.
Then
\be\label{tmp}
T(t,r)=q^N,\quad T(0,r-1)=S,\quad T(0,r)=S^*,
\ee
and $T(j,i)=I$ for all other pairs $(j,i)$, except possibly
$T(t-1,r-1)$,  which is $S^*$
if $t-1\in J_{r-1}$, and $I$ otherwise.
It is easy to check that for $\gamma\in \scrf_{\gamma_i}$,
$\gamma(t,r)+\gamma(0,r)=\gamma_i(t,r)+\gamma_i(0,r)$.
Therefore the set
$\{\gamma(t,r)+\gamma(0,r):\gamma\in\cup_i\scrf_{\gamma_i}\}$
is bounded. Let $n\in\bbn$ be such that this set is contained in $[-n,n]$.
Suppose $\gamma\in\cup_i\scrf_{\gamma_i}$ obey $\gamma(t,r)=0$.
Then it follows from (\ref{tmp}) that
$T^{2n+1}e_\gamma=e_{\gamma'}$, where
\[
 \gamma'(0,r)=\gamma(0,r)+2n+1,\quad
   \gamma'(0,r-1)=\gamma(0,r-1)-2n-1,\quad
    \gamma'(t,r)=\gamma(t,r).
\]
It is clear from this that $\gamma'\not\in\cup_i\scrf_{\gamma_i}$,
so that $PT^{2n+1}e_\gamma=0$. This means
$[P,T^{2n+1}]e_\gamma=-e_{\gamma'}$ for all $\gamma\in\cup_i\scrf_{\gamma_i}$
with $\gamma(t,r)=0$. Since there are infinitely many choices of such $\gamma$,
it follows that  $[P,T^{2n+1}]$ can not be compact.

We thus have the following theorem.
\bthm
Let $\ell>1$. Then there does not exist any
Dirac operator on $L_2(\Gamma)$ that diagonalises
with respect to the canonical orthonormal basis
and has nontrivial sign.
\ethm

\brmrk
{}Let $F$ be a subset of $\{1,2,\ldots,\ell\}$. Define
$\pi_{\omega,F}$ to be the representation obtained by integrating
$\psi_\omega\ast \rchi_\bldz$ with respect to those components $z_i$
of $\bldz$ for which $i\in F$. If one looks at the representations
$\pi_{\omega,F}$ instead of $\pi_\omega$,
a similar analysis will show that nontrivial spectral triples would exist
only in the case  where $\omega$ is of the form $s_k$ (so that $\ell(\omega)=1$),
and $F=\{k\}$. The nontrivial triples in this case will essentially be those of $SU_q(2)$
obtained in \cite{c-p2} and will correspond to the `$k$\raisebox{.4ex}{th} copy'
of $SU_q(2)$ sitting inside $SU_q(\ell+1)$ via
the map
\[
u_{ij}\mapsto \cases{ \alpha & if $j=i=k$,\cr
                      \alpha^* & if $j=i=k+1$,\cr
                      -q\beta^* & if $j=k+1, i=k$,\cr
                       \beta & if $j=k, i=k+1$,\cr
                       I & if $i=j$,\cr
                       0 & otherwise.}
\]
\ermrk

{\footnotesize \textbf{Acknowledgement.} We would like to thank
Prof.\ Yan Soibelman for his remarks on an earlier paper
on $SU_q(2)$ which encouraged us to look at the
$SU_q(\ell+1)$ case.
}


\noindent{\sc Partha Sarathi Chakraborty}
(\texttt{chakrabortyps@cf.ac.uk})\\
         {\footnotesize School of Mathematics,
 Cardiff University, Senghennydd Road, Cardiff, UK}\\[1ex]
{\sc Arupkumar Pal} (\texttt{arup@isid.ac.in})\\
         {\footnotesize Indian Statistical
Institute, 7, SJSS Marg, New Delhi--110\,016, INDIA}


\end{document}